\input amstex 
\documentstyle{amsppt} 
\input xy
\xyoption{all}
\loadbold
\magnification=1200
\NoBlackBoxes
\pagewidth{6.4truein}
\pageheight{8.51truein}

\redefine\emptyset{{\varnothing}}
\define\ac{\acuteaccent}

\redefine\cdot{\boldsymbol\cdot}
\define\Pre{\operatorname{Pre}}
\define\Shv{\operatorname{Shv}}

\NoRunningHeads

\topmatter 

\title  Excision for simplicial sheaves on the Stein site 
 \\ and Gromov's Oka Principle   \endtitle
\author Finnur L\ac arusson  \endauthor
\affil University of Western Ontario \endaffil
\address Department of Mathematics, University of Western Ontario,
        London, Ontario N6A~5B7, Canada \endaddress
\email larusson\@uwo.ca \endemail

\thanks The author was supported in part by the Natural Sciences and 
Engineering Research Council of Canada. \endthanks

\thanks Completed December 2000; minor changes April 2002. \endthanks

\abstract A complex manifold $X$ is said to satisfy the Oka-Grauert 
property if the inclusion $\Cal O(S,X)\hookrightarrow\Cal C(S,X)$ is 
a weak equivalence for every Stein manifold $S$, where the spaces of
holomorphic and continuous maps from $S$ to $X$ are given the
compact-open topology.  Gromov's Oka principle states that if $X$ has a
spray, then it has the Oka-Grauert property.  The purpose of this paper
is to investigate the Oka-Grauert property using homotopical algebra. 
We embed the category of complex manifolds into the model category of
simplicial sheaves on the site of Stein manifolds.  Our main result is
that the Oka-Grauert property is equivalent to $X$ representing a finite
homotopy sheaf on the Stein site.  This expresses the Oka-Grauert
property in purely holomorphic terms, without reference to continuous
maps.  \endabstract


\endtopmatter

\document

\specialhead 1. Introduction \endspecialhead

\noindent The purpose of this paper is to investigate Gromov's Oka
principle using abstract homotopy theory and to recast it in intrinsic,
holomorphic terms, and thereby to introduce some of the concepts and
methods of homotopical algebra into complex geometry. 

The Oka principle is a vague maxim, supported by many results.  It may
be phrased by saying that on a Stein manifold (complex submanifold of
Euclidean space), analytic problems of a cohomological nature have only
topological obstructions.  It has a long and venerable history, starting
with the 1939 result of Oka stating, in modern language, that a
holomorphic line bundle on a domain of holomorphy is trivial if it is
topologically trivial.  Deep generalizations to vector bundles and
certain other fiber bundles were obtained by Grauert in the late 1950s. 
Another manifestation of the Oka principle is the 1975 result of
Cornalba and Griffiths that every rational cohomology class of degree
$2k$ on a Stein manifold is a rational multiple of the fundamental class
of a $k$-codimensional analytic subvariety.  For a survey, see
\cite{19}. 

A major development appeared in Gromov's 1989 paper \cite{13}.  He
discovered that if a complex manifold $X$ has a geometric structure
called a {\it spray}, then $X$ satisfies what we shall call the {\it 
Oka-Grauert property}, meaning that the inclusion $\Cal O(S,X)
\hookrightarrow\Cal C(S,X)$ is a weak equivalence for every Stein 
manifold $S$, where the spaces of holomorphic and continuous maps from 
$S$ to $X$ are given the compact-open topology.  In particular, 
every continuous map from a
Stein manifold to $X$ can be deformed to a holomorphic map.  A spray on
$X$ consists of holomorphic maps $\Bbb C^m\to X$, $0\mapsto x$, one for
each $x\in X$, submersive at $0$, and varying holomorphically with $x$. 
For a detailed proof of Gromov's theorem and a thorough introduction,
see \cite{8}.  A more general version for sections of submersions is
contained in \cite{10}.  Gromov's Oka principle has been applied to the
famous problem of embedding Stein manifolds into Euclidean spaces of the
smallest possible dimension: Sch\"urmann has used it to prove Forster's
conjecture in higher dimensions \cite{28}, following work of Eliashberg
and Gromov \cite{3}.  Other applications and variants may be found in 
\cite{4, 5, 6, 7, 9}. 

The Oka-Grauert property certainly has a homotopy-theoretic flavour. 
Our goal is to turn this impression into a precise statement in an
abstract homotopy-theoretic setting.  At the same time, we will express
the Oka-Grauert property in purely holomorphic terms, without reference
to continuous maps.

\smallskip Abstract homotopy theory, also known as homotopical algebra,
is due to Quillen \cite{2, 12, 14, 16, 27}.  Its fundamental notion
is that of a {\it model category}: a category satisfying certain axioms
that allow us to develop an analogue of ordinary homotopy theory.  There
has been much activity in recent years in both the theory and applications
of homotopical algebra, the most notable example being the development of 
motivic homotopy theory, i.e., the homotopy theory of schemes, leading to 
Voevodsky's proof of the Milnor conjecture \cite{31, 32}.  As far as I 
know, the present paper is the first attempt to introduce homotopical 
algebra into analytic geometry. 

The first step, just as for schemes, is to embed the category of complex
manifolds into the model category of simplicial objects in a topos by a
Yoneda embedding of some sort, where we can then do homotopy theory with
them.  Here, this is done in Section 4.  Whereas in algebraic geometry
the focus is on generalized cohomology theories and ultimately motives,
on the analytic side it seems of more immediate interest to try to do
algebraic topology with complex manifolds and holomorphic maps instead
of topological spaces and continuous maps, and then our attention is
immediately drawn to Gromov's Oka principle.  For the purpose of
studying the Oka principle, we associate to a complex manifold $X$ the
simplicial sheaf $s\Cal O(\cdot,X)$ on the site of all Stein manifolds. 
Here, spaces of holomorphic maps are given the compact-open topology; it
is for technical reasons that we turn them into simplicial sets by
applying the singular functor $s$.  Using the compact-open topology
allows us to work at a relatively simple technical level: we neither
have to localize nor stabilize to get something interesting.  

Among other things, a new notion of weak equivalence between complex 
manifolds emerges, weaker than biholomorphism and stronger than 
topological weak equivalence (it can tell apart the punctured plane and 
the punctured disc).  The definition is simple --- a holomorphic map 
$X\to Y$ is a weak equivalence if it induces a topological weak 
equivalence $\Cal O(S,X)\to\Cal O(S,Y)$ for every Stein manifold $S$ 
--- but the point is that it fits into a model structure.

The main result of this paper is Theorem 2.1, later rephrased as Theorem
4.3.  It states that a complex manifold $X$ has
the Oka-Grauert property if and only if the
simplicial sheaf $s\Cal O(\cdot,X)$ is a {\it finite homotopy sheaf}
on the Stein site.  This homotopy-theoretic property is also called {\it
finite excision}.  It gives rise to Mayer-Vietoris sequences of homotopy 
groups and is familiar from topology and appears nowadays in algebraic
geometry: see e.g\.  \cite{25, \S 3.1.2}.  I have tried to make the
proof of Theorem 2.1 as understandable as possible to those unfamiliar
with homotopy theory.  However, the very definition of excision requires
the notion of a homotopy limit (a deformation invariant approximate
limit), for which I refer the reader to \cite{12, VIII.2} and \cite{14,
Ch\.  19}.  The proof uses the main result of Section 3, Theorem 3.4,
whose crucial ingredient is a classical theorem of Brown and Gersten
\cite{1, Thm\.  1}, foundational in the homotopy theory of simplicial
sheaves.  Section 3 is pure homotopy theory and constitutes the bulk of
the paper.  The proof of Theorem 2.1 also uses Siu's theorem on the
existence of Stein neighbourhoods of Stein subvarieties, and Whitney's
lemma on decomposing an open subset in Euclidean space into a union of
cubes with special properties.  

Those familiar with the homotopy theory of simplicial sheaves would now
ask if the Oka-Grauert property is actually equivalent to descent,
making Gromov's Oka principle somewhat analogous to such results as
Brown-Gersten, Nisnevich, and Thomason descent in algebraic geometry
\cite{24}.  Descent may be loosely described as a homotopic 
local-to-global property, stronger than excision.  I do not know the 
answer: the finiteness properties that make descent possible in algebra
do not hold in analysis.  I hope to address this question in future work.
In the meantime, Section 5 contains a partial descent theorem of sorts 
for quasi-projective manifolds. 

It should be emphasized that this paper is not about the proof of
Gromov's theorem at all.  The ideas presented in the paper may shed new 
light on the proof, but this possibility is not pursued here.  Rather, 
the paper is about the Oka-Grauert property itself and its 
homotopy-theoretic meaning in a new model-categorical context 
for complex manifolds.

A word about terminology: we take manifolds to be second countable by
definition, but not necessarily connected.

\smallskip
\noindent {\it Acknowledgements.} I am grateful for discussions with
Paul Balmer, Dan Christensen, Gaunce Lewis, Sergei Yagunov, and
especially Rick Jardine, who has generously and patiently answered many
questions of mine.

\vfill

\specialhead 2. The Oka-Grauert property is equivalent to finite excision
\endspecialhead

\noindent 
We say that a complex manifold $X$ (second countable but not necessarily
connected) has the {\it Oka-Grauert property} if the inclusion map
$$\Cal O(S,X)\hookrightarrow \Cal C(S,X)$$
is a weak equivalence for all Stein manifolds $S$, where the spaces of
holomorphic and continuous maps from $S$ to $X$ carry the
compact-open topology.  This means that the inclusion induces
isomorphisms of all homotopy groups for all base points, as well as a
bijection of path components.  In particular, surjectivity on path
components means that every continuous map $S\to X$ is homotopic to a
holomorphic map.  Note that requiring $S$ to be connected results in an 
equivalent condition.

We say that $X$ satisfies {\it finite excision} if
whenever $\{U_1,\dots,U_m\}$ is a finite cover of a Stein manifold $S$
by Stein open subsets, $\Cal O(S,X)$ is not only the limit but also the
homotopy limit of the diagram whose objects are the spaces $\Cal
O(U_{i_1}\cap\dots\cap U_{i_k}, X)$ for $1 \leq i_1 < \dots < i_k \leq
m$ and $k=1,\dots,m$, and whose arrows are the restriction maps.  This
diagram forms an $m$-dimensional cube with one vertex missing; the limit
or homotopy limit provides the missing vertex.  For $m=2$, this simply
means that the square
$$\CD
\Cal O(S,X) @>>> \Cal O(U_1,X)  \\
@VVV               @VVV           \\
\Cal O(U_2,X) @>>> \Cal O(U_1\cap U_2, X)
\endCD $$
is not only a pullback but also a homotopy pullback.

Let us make clear what we mean by a homotopy limit.  We view homotopy
limits as determined only up to weak equivalence.  We say that a
topological space $Y$ with a map to a diagram $\Cal X$ of spaces over an
index category $J$ is the homotopy limit of $\Cal X$ (or more properly
that the diagram $Y\to\Cal X$ is a homotopy limit) if and only if the
singular set $sY$ of $Y$ is the homotopy limit of the diagram $s\Cal X$
of simplicial sets.  This in turn means that if $\Cal A$ is a fibrant
model for $s\Cal X$ in the category of diagrams of simplicial sets over
$J$ with the pointwise cofibration structure, then the natural map from
$sY$ to the limit of $\Cal A$ is a weak equivalence \cite{12,
VIII.2.11}.  It is my understanding that this definition is standard (up
to weak equivalence): it agrees up to weak equivalence with the
definitions of Bousfield-Kan and Hirschhorn \cite{14, Ch\.  19}.  The
homotopy limit of a diagram of spaces or simplicial sets may also be
described somewhat explicitly as the function space or complex of
morphisms to the diagram from a certain diagram associated to the index
category \cite{12, VIII.2.3; 14, 19.1.10}.  We will make frequent use
of the basic fact that if one of two weakly equivalent diagrams of
spaces or fibrant simplicial sets is a homotopy limit, then so is the
other one. 

\smallskip   
Our main result is that the Oka-Grauert property is equivalent to 
finite excision. 

\proclaim{2.1. Theorem}  A complex manifold has the Oka-Grauert
property if and only if it satisfies finite excision.
\endproclaim

\demo{Proof} Let $X$ be a complex manifold with the Oka-Grauert
property.  Consider a finite cover of a Stein manifold $S$ by Stein open
subsets $U_1,\dots,U_m$.  We have two $m$-dimensional cube diagrams of
spaces, one with objects $\Cal O(U_{i_1}\cap\dots\cap U_{i_k}, X)$, and
the other with objects $\Cal C(U_{i_1}\cap\dots\cap U_{i_k}, X)$ for $1
\leq i_1 < \dots < i_k \leq m$.  We have a morphism of inclusions from
the first diagram to the second one, consisting of weak equivalences by
assumption (here we need to know that the intersection of Stein open
sets is Stein; see Lemma 4.1).  Hence, the induced map between the
homotopy limits of the two diagrams is a weak equivalence \cite{12,
VIII.2.2}.  By Theorem 3.4 below, the homotopy limit of the second
diagram is $\Cal C(S, X)$, which by assumption is weakly equivalent to
$\Cal O(S,X)$. 

Conversely, assume $X$ is a complex manifold satisfying finite excision and 
let $S$ be a Stein manifold.  We first reduce our problem to the case when
$S$ is a domain in Euclidean space.  As noted above, we may take $S$ to
be connected, so $S$ embeds into some Euclidean space.  Then, by a
theorem of Siu \cite{29, Cor\.  1}, there is a connected Stein
neighbourhood $V$ of $S$ and a holomorphic retraction $\rho:V\to S$. 
Let $\iota:S\hookrightarrow V$ be the inclusion, so $\rho\circ\iota =
\text{id}_S$, and we have a diagram
$$\xymatrix{
\Cal O(S,X) \ar@<1ex>[r]^{\rho^*} \ar[d]_\phi & 
\Cal O(V,X) \ar@<1ex>[l]^{\iota^*} \ar[d]^\psi \\
\Cal C(S,X) \ar@<1ex>[r]^{\rho^*} & \Cal C(V,X) \ar@<1ex>[l]^{\iota^*}  } $$
where $\phi$ and $\psi$ are the inclusions.  Now suppose $\psi$ is a
weak equivalence.  Since $\iota^*\circ\rho^*=\text{id}$, $\rho^*$
induces monomorphisms and $\iota^*$ induces epimorphisms on all homotopy
groups.  Hence, $\rho^*\circ\phi=\psi\circ\rho^*$ induces monomorphisms
on all homotopy groups so $\phi$ does too, and 
$\phi\circ\iota^*=\iota^*\circ\psi$ induces epimorphisms on all homotopy
groups so $\phi$ does too, and $\phi$ is a weak equivalence.

To complete the proof we need to show that $\Cal O(V,X)\hookrightarrow
\Cal C(V,X)$ is a weak equivalence when $V$ is a Stein domain in
Euclidean space.  We want to express $V$ as a finite union of open
subsets all of whose connected components are convex.  This can surely
be done in many ways.  We shall refer to Whitney's classical lemma on
decomposing an open set in Euclidean space into a union of cubes with
special properties \cite{30, VI.1}.  We get that $V=U_1\cup\dots\cup
U_m$, where each $U_i$ is a disjoint union of open cubes with sides
parallel to the coordinate axes.  (This is not stated explicitly in
\cite{30}, but may easily be obtained from there.) Then every
intersection $U=U_{i_1}\cap\dots\cap U_{i_k}$ is a disjoint union of
open boxes and hence Stein.  Each box is {\it holomorphically
contractible} in the sense that the identity map can be joined to a
constant map by a continuous family of holomorphic maps, so the
inclusion $\Cal O(U,X)\hookrightarrow \Cal C(U,X)$ is clearly a weak
equivalence.  Now we look at two $m$-dimensional cube diagrams of
spaces, one with objects $\Cal O(U_{i_1}\cap\dots\cap U_{i_k}, X)$ and
the other with objects $\Cal C(U_{i_1}\cap\dots\cap U_{i_k}, X)$ for $1
\leq i_1 < \dots < i_k\leq m$.  As above, we see that the induced map
between the homotopy limits of the two diagrams is a weak equivalence. 
Since $X$ satisfies finite excision, and by Theorem 3.4, this map is the
inclusion $\Cal O(V,X)\hookrightarrow \Cal C(V,X)$ (at least up to weak
equivalence), and the proof is complete.  \qed\enddemo

The proof shows that a complex manifold $X$ satisfies finite excision
if and only if $\Cal O(S,X)\to \Cal O(\Cal B,X)$ is a homotopy 
limit for every Stein basis $\Cal B$ for a binoetherian subtopology on a 
Stein manifold $S$, viewed as a subdiagram of the site of $S$. 

The proof also shows that the Oka-Grauert property for a complex
manifold $X$ is equivalent to the inclusion $\Cal
O(V,X)\hookrightarrow\Cal C(V,X)$ being a weak equivalence for all
domains of holomorphy $V$ in $\Bbb C^n$ for all $n\geq 1$.  To some
extent it is therefore a matter of taste whether one chooses to work
with Stein manifolds or domains of holomorphy in Euclidean space in the
present context. 

\smallskip
If $X$ has the Oka-Grauert property and $S$ is a Stein manifold, it is
natural to ask whether the weak equivalence $\Cal O(S,X)\hookrightarrow
\Cal C(S,X)$ is actually a homotopy equivalence.  When $S$ is algebraic,
a little topology shows that the inclusion has a right homotopy inverse,
so there is a continuous way of associating to each continuous map $S\to
X$ a holomorphic map homotopic to it.  Note that we do not assert that
the homotopy inverse fixes holomorphic maps. 

\proclaim{2.2.  Theorem} Let $S$ be an affine algebraic manifold, i.e.,
a Stein manifold biholomorphic to an algebraic submanifold of Euclidean
space.  If $X$ is a complex manifold and the inclusion $\Cal O(S,X)
\hookrightarrow \Cal C(S,X)$ is a weak equivalence, then it has a right
homotopy inverse.  
\endproclaim

\demo{Proof} Being a smooth manifold, $X$ has a countable triangulation
\cite{26, 10.6}, so $X$ is homeomorphic to a countable CW complex. 
Also, $S$ is homotopy equivalent to a finite CW complex $K$ \cite{23,
Lemma A.3}, and $\Cal C(S,X)$ is homotopy equivalent to $\Cal C(K,X)$,
which has the homotopy type of a (countable) CW complex \cite{11,
5.2.5}, so $\Cal C(S,X)$ has the homotopy type of a CW complex.  Hence,
the natural map $a:|s\Cal O(S,X)|\to |s\Cal C(S,X)|\to \Cal C(S,X)$,
which is a weak equivalence by assumption, has a homotopy inverse $b$. 
Let $c$ be the natural map $|s\Cal O(S,X)|\to \Cal O(S,X)$.  By
adjunction, $a=ic$, where $i$ is the inclusion $\Cal O(S,X)
\hookrightarrow \Cal C(S,X)$.  Now $i(cb)=ab$ is homotopic to the
identity on $\Cal C(S,X)$, so $cb$ is a right homotopy inverse for $i$. 
\qed\enddemo

\specialhead 3.  Excision and Brown-Gersten descent \endspecialhead

\noindent The main purpose of this section is to establish the excision
property of $s\Cal C(\cdot,X)$ used in the previous section (Theorem
3.4).  One feels that it should be possible to verify this directly
using the explicit description of the homotopy limit given in \cite{12,
VIII.2.3} and \cite{14, 19.1.10}, but rather than attempt this, we give
a proof based on Brown-Gersten descent (Theorem 3.1).  Brown-Gersten
descent is surely well known among experts, but in the absence of a good
reference, we have provided a detailed proof.  Theorems 3.3 and 3.4 are
new as far as I know.  We start with a brief review of the basic notions
of the homotopy theory of simplicial presheaves. 

Let $\Cal S$ be a small Grothendieck site.  A {\it simplicial presheaf}
on $\Cal S$ is a contravariant functor from $\Cal S$ to the category
$s\bold{Set}$ of simplicial sets.  There is a standard model structure
on the category $s\Pre\Cal S$ of simplicial presheaves on $\Cal S$ in
which the cofibrations are monomorphisms, i.e., pointwise injections
(where {\it pointwise} means {\it at every object of the site}), and a
weak equivalence is a map that induces isomorphisms of all homotopy
sheaves \cite{17, p\.  59}.  If $\Cal S$ has enough points, e.g\.  if
it is the site of a topological space, then this is equivalent to the
map inducing weak equivalences of all stalks.  A weak equivalence is
still a weak equivalence with respect to any finer topology.  Fibrations
are defined by a right lifting property.  There is an induced model
structure on the full subcategory $s\Shv\Cal S$ of simplicial sheaves on
$\Cal S$.  More is true: both categories are {\it proper, simplicial
model categories} \cite{17, 18}.  There is another model structure
on $s\Pre\Cal S$ given by the trivial topology on $\Cal S$, in which the
only covers are those consisting of a single isomorphism.  In this
structure, the weak equivalences are the pointwise weak equivalences. 
The words {\it fine} and {\it finely} shall refer to the former model
structure, and the words {\it coarse} and {\it coarsely} to the latter. 
A map of simplicial presheaves which is a weak equivalence will be
referred to as {\it acyclic}.  (We are trying not to overuse the word
{\it trivial}.) A coarse weak equivalence is a fine weak equivalence; a
fine fibration is a coarse fibration. 

The concept of a fibrant object is fundamental in homotopy theory.  We
will now define several important weaker notions and briefly describe
their relationships. 

We say that a simplicial presheaf $G$ on $\Cal S$ satisfies {\it
descent} if any fine weak equivalence from $G$ to a finely fibrant
simplicial presheaf on $\Cal S$ is coarsely acyclic.  Equivalently
(using the Whitehead Theorem that a weak equivalence between bifibrant
objects is a homotopy equivalence), a finely fibrant model for $G$ is
also a coarsely fibrant model for $G$.  This notion is invariant under
coarse weak equivalences.  A finely fibrant simplicial presheaf
satisfies descent. 

It may be shown that a coarsely fibrant simplicial presheaf $G$ on $\Cal
S$ is both pointwise fibrant and {\it flabby} (or {\it flasque}), which
means that the restriction map $G(U)\to G(V)$ is a fibration for every
monomorphism $V\to U$ in $\Cal S$. 

Now let $X$ be a topological space.  We say that a pointwise fibrant
simplicial presheaf $G$ on $X$, i.e., a presheaf of Kan complexes,
satisfies {\it excision} (or {\it two-set excision}) if $G(\emptyset)$
is contractible (this is true if $G$ is a sheaf) and whenever $U$ and
$V$ are open in $X$, the square
$$ \CD
G(U\cup V) @>>> G(U) \\
@VVV            @VVV \\
G(V)   @>>>  G(U\cap V)
\endCD $$
is a homotopy pullback.  This notion was introduced by Brown and Gersten
\cite{1, \S 2}, who used the term {\it pseudo-flasque}.  It is
clearly invariant under coarse weak equivalences.  If $G$ is a flabby
simplicial sheaf (so $G$ is in particular pointwise fibrant), then $G$
satisfies excision: the diagram $G(U)\to G(U\cap V)\gets G(V)$ is
fibrant in its diagram category, so the ordinary pullback $G(U\cup V)$
is the homotopy pullback \cite{12, VI.1.8}.  If $G$ is a pointwise
fibrant simplicial presheaf satisfying descent, consider the fine weak
equivalence from $G$ to its sheafification $aG$, and the fine weak
equivalence from $aG$ to a finely fibrant model $F$ in $s\Shv X$.  Then
$F$ is also finely fibrant in $s\Pre X$, so $F$ is flabby and satisfies
excision, and since $F$ and $G$ are finely and hence coarsely
equivalent, $G$ does too.  Under a strong finiteness condition on $X$,
the converse holds: this is Brown-Gersten descent.  We say that a
topological space is {\it binoetherian} if both the open sets and the
irreducible closed sets satisfy the ascending chain condition.  An
important example is a space with a finite topology. 

\proclaim{3.1.  Theorem (Brown-Gersten descent)} For pointwise fibrant
simplicial presheaves on a binoetherian space, excision is equivalent to
descent.  
\endproclaim

\demo{Proof} Let $G$ be a pointwise fibrant simplicial presheaf on a
binoetherian space, and suppose that $G$ satisfies excision.  Let $G\to
F$ be a fine weak equivalence from $G$ to a finely fibrant simplicial
presheaf $F$.  Then $F$ satisfies descent and hence excision.  The
theorem now follows from the next result.  
\qed\enddemo

\proclaim{3.2.  Proposition} Let $F$ and $G$ be pointwise fibrant
simplicial presheaves satisfying excision on a binoetherian space $X$. 
Then a fine weak equivalence $G\to F$ is coarsely acyclic.  
\endproclaim

Our argument is an adaptation and explication of Morel and Voevodsky's
proof of unstable Nisnevich descent in \cite{25, \S 3.1.2}. 

\demo{Proof} Let $F'$, $G'$ be coarsely fibrant models in $s\Pre X$ for
$F$, $G$ respectively.  Factor the induced map $G'\to F'$ as a coarse
weak equivalence $G'\to G''$ followed by a coarse fibration $G''\to F'$. 
Then $G''$ is coarsely fibrant and it suffices to show that the map
$G''\to F'$ is coarsely acyclic.  By replacing $F$, $G$ by $F'$, $G''$,
we may assume that $F$ and $G$ are coarsely fibrant and hence flabby and
pointwise fibrant, and that the fine weak equivalence $G\to F$ is a
coarse fibration and hence a pointwise fibration. 

It suffices to show that for any open set $U$ in $X$ and any vertex $x$
in $F(U)$, the fibre of the fibration $G(U)\to F(U)$ over $x$ is
contractible (in particular nonempty): then $G(U)\to F(U)$ is a weak
equivalence (if $F(U)$ is empty, then so is $G(U)$ and this is still
true).  Note that $U$ is binoetherian in the subspace topology, so we
may assume that $U=X$.  Fix a vertex $x$ in $F(X)$ and consider the
simplicial presheaf $K$ on $X$ that associates to an open set $V$ in $X$
the fibre of $G(V)\to F(V)$ over the image of $x$ in $F(V)$.  Note that
$K$ is a pullback of a diagram $G\to F\gets\ast$ of simplicial
presheaves on $X$, where $\ast$ denotes the final simplicial presheaf. 
Hence, $K$ is coarsely fibrant, and therefore flabby and pointwise
fibrant. 

We need to show that $K$ is pointwise contractible.  Since $K$ is
stalkwise contractible (pullbacks commute with filtered colimits, so
taking fibres commutes with taking stalks), this follows from \cite{1,
Thm\.  1} along with the remark at the end of \cite{1, \S 2}, once we
know that $K$ satisfies excision.  To complete the proof, let us verify
this. 

First of all, since $F(\emptyset)$ and $G(\emptyset)$ are contractible,
so is the fibre $K(\emptyset)$.  Now let $U$ and $V$ be open in $X$. 
Since $K$ is flabby and pointwise fibrant, the homotopy pullback of
$K(U)\to K(U\cap V)\gets K(V)$ is the ordinary pullback, so we need to
show that the natural map $K(U\cup V)\to K(U)\times_{K(U\cap V)} K(V)$
is a weak equivalence.  By the cogluing lemma \cite{12, II.8.13} applied
to the natural map from the pullback square
$$\CD
K(U\cup V) @>>> G(U\cup V)  \\
@VVV          @VVV          \\
\ast  @>>>  F(U\cup V)
\endCD $$
to the pullback square
$$\CD
K(U)\times_{K(U\cap V)} K(V)  @>>>  G(U)\times_{G(U\cap V)} G(V)  \\
@VVV					@VpVV                      \\
\ast @>>> F(U)\times_{F(U\cap V)} F(V)
\endCD $$
we would be done if $p$ was a fibration.  Since this is not to be
expected, we need to replace $p$ by a fibration in a reasonable way. 

We shall work in the category of squares of simplicial sets of the type
$$\CD
4 @>>> 3 \\
@VVV @VVV \\
2 @>>> 1
\endCD $$
This category carries the pointwise cofibration and the pointwise
fibration simplicial model structures \cite{12, p\. 403}.  In both of 
them, the weak equivalences are the pointwise weak equivalences.  We 
shall write $Q_F$ for the square
$$\CD
F(U\cup V) @>>> F(V)  \\
@VVV	      @VVV    \\
F(U)  @>>>  F(U\cap V)
\endCD $$
and similarly for $G$ and $K$.  In the pointwise cofibration structure,
factor $Q_G \to Q_F$ as a weak equivalence $Q_G \to Q$ followed by a
fibration $Q\to Q_F$.  Using the fibration $Q\to Q_F$, one can show that
$Q$ is pointwise fibrant and the maps $Q_2, Q_3\to Q_1$ are fibrations. 
Since $Q_G$ is a homotopy pullback, so is $Q$, and the natural map $Q_4
\to Q_2\times_{Q_1}Q_3$ is a weak equivalence. 

The cogluing lemma applied to the natural map from the pullback square
$$\CD
Q_K @>>> Q_G  \\
@VVV  @VVV    \\
\ast @>>> Q_F
\endCD $$
to the pullback square
$$\CD
P @>>> Q  \\
@VVV  @VVV  \\
\ast @>>> Q_F
\endCD $$
with respect to the pointwise fibration structure (which is right
proper) shows that $Q_K\to P$ is a weak equivalence (we invoke the fact
that a fibration in the pointwise cofibration structure is a pointwise
fibration).  Also, $P$ is fibrant in the pointwise cofibration
structure.  Hence, it suffices to show that the square $P$ is a homotopy
pullback, which is the case if and only if the natural map $P_4\to
P_2\times_{P_1}P_3$ is a weak equivalence. 

By the cogluing lemma applied to the natural map from the pullback square
$$\CD
P_4 @>>> Q_4  \\
@VVV          @VVV          \\
\ast  @>>>  F(U\cup V)
\endCD $$
to the pullback square
$$\CD
P_2\times_{P_1} P_3  @>>>  Q_2\times_{Q_1} Q_3  \\
@VVV					@VpVV       \\
\ast @>>> F(U)\times_{F(U\cap V)} F(V)
\endCD $$
it suffices to show that the map $p$ is a fibration.  Therefore, to
complete the proof, we need to show that if $R\to S$ is a fibration 
of squares in the pointwise cofibration structure, then 
$R_2\times_{R_1} R_3 \to S_2\times_{S_1} S_3$ is a fibration of 
simplicial sets.  (Having a pointwise fibration is not enough.)  
A square
$$\CD
A @>>> R_2\times_{R_1} R_3 \\
@VVV          @VVV         \\
B @>>> S_2\times_{S_1} S_3
\endCD $$
is the same thing as a square
$$\CD
\tilde A @>>> R \\
@VVV  @VVV      \\
\tilde B @>>> S
\endCD $$
of squares, where $\tilde A$ is the square
$$\xymatrix{
\emptyset \ar[r] \ar[d] & A \ar@{=}[d] \\ 
A \ar@{=}[r] & A
} $$
and $\tilde B$ is defined similarly.  A map $A\to B$ is an acyclic
cofibration if and only if the induced map $\tilde A\to\tilde B$ is a
pointwise acyclic cofibration, and then a lifting in the latter square
gives a lifting in the former.
\qed\enddemo

Using Brown-Gersten descent, we can now strengthen the excision condition.

\proclaim{3.3.  Theorem} Let $G$ be a pointwise fibrant simplicial
presheaf satisfying excision on a topological space $X$.  Let $\Cal B$
be a basis for a binoetherian subtopology on an open set $U$ in $X$,
viewed as a subdiagram of the site of $X$.  Then $G(U)$ is the homotopy
limit of $G|\Cal B$.  
\endproclaim

Recall that a collection $\Cal B$ of open subsets in a topological space
is a basis for a subtopology on $\bigcup\Cal B$ if whenever $U,V\in\Cal
B$ and $p\in U\cap V$, there is $W\in\Cal B$ with $p\in W\subset U\cap
V$.  We do not assume that a basis is closed under intersections.  The
definition of excision refers to the case when $\Cal B$ consists of two
open sets in $X$ and their intersection. 

\demo{Proof} We may assume that $X=U$ is binoetherian and $\Cal B$ is a
basis for the topology on $X$.  By Brown-Gersten descent, $G$ satisfies
descent, so there is a coarse weak equivalence from $G$ to a finely
fibrant simplicial {\it sheaf} $F$ on $X$.  Let us show that $F|\Cal B$
is fibrant in the diagram category $s\bold{Set}^{\Cal B^{\text{op}}}$
with the pointwise cofibration structure.  Let $A\to B$ be an acyclic
cofibration in $s\bold{Set}^{\Cal B^{\text{op}}}$ and $A\to F|\Cal B$ be
a morphism.  Now $A$ and $B$ yield \ac etale spaces over $X$ with
sheaves of sections $\tilde A$ and $\tilde B$ respectively, such that
the induced map $\tilde A\to\tilde B$ is a finely acyclic cofibration in
$s\Pre X$.  Of course the diagram $F|\Cal B$ yields the sheaf $F$ itself
in this way, so the map $A\to F|\Cal B$ factors through $\tilde A|\Cal
B$.  Since $F$ is finely fibrant, $\tilde A\to F$ factors through
$\tilde B$, so $A\to F|\Cal B$ factors through $B$.  This shows that
$F|\Cal B$ is a fibrant model for $G|\Cal B$ in $s\bold{Set}^{\Cal
B^{\text{op}}}$.  Hence, the homotopy limit of $G|\Cal B$ is weakly
equivalent to the limit of $F|\Cal B$, which is $F(X)$ since $F$ is a
sheaf, and $F(X)$ is weakly equivalent to $G(X)$.  \qed\enddemo

The following theorem gives the excision property of $\Cal C(\cdot,X)$
used in the previous section.  We shall make brief use of the category
$\bold{Space}$ of compactly generated weak Hausdorff spaces \cite{11,
A.1; 22, Ch\.  5}, \lq\lq the category of spaces in which algebraic
topologists customarily work\rq\rq\ \cite{22, p\.  37}.  We denote the
internal function complex in $s\bold{Set}$ by
$\bold{Hom}_{s\bold{Set}}(\cdot,\cdot)$ as in \cite{12, I.5}. 

\proclaim{3.4.  Theorem} Let $X$ be a smooth manifold and $Y$ be a
compactly generated weak Hausdorff space.  If $\Cal B$ is a basis for a
binoetherian subtopology on an open set $U$ in $X$, viewed as a
subdiagram of the site of $X$, then $\Cal C(U,Y)$ is the homotopy limit
of $\Cal C(\cdot,Y)|\Cal B$. 
\endproclaim

Here, {\it smooth} means at least once continuously differentiable.  For
this result, the most important consequence of $X$ being a smooth
manifold is that every open subset of $X$ is cofibrant since it has a
triangulation \cite{26, 10.6}.  Surely, the class of spaces with this
property is much larger than the class of smooth manifolds, but I am not
aware of any description of it.  We also need every open subset of $X$
to be normal. 

\demo{Proof} By Theorem 3.3, we need to verify that the simplicial sheaf
$s\Cal C(\cdot,Y)$ on $X$ satisfies excision.  We will show that the
simplicial presheaf $F=\bold{Hom}_{s\bold{Set}}(s\cdot,sY)$ on $X$ is
pointwise fibrant and satisfies excision and that $F$ is coarsely weakly
equivalent to $s\Cal C(\cdot,Y)$.  Note first that $F(\emptyset)$ is the
singular set of a point, so $F(\emptyset)$ is contractible. 

Consider the functor $H=\bold{Hom}_{s\bold{Set}}(\cdot,sY)$ from
$s\bold{Set}^{\text{op}}$ to $s\bold{Set}$.  By Quillen's Axiom SM7 for
a simplicial model category \cite{12, II.3.1}, if $A\to B$ is a
cofibration in $s\bold{Set}$, then the induced map $H(B)\to H(A)$ is a
fibration.  Hence, $F$ is flabby and pointwise fibrant.  If we knew that
$F$ was a sheaf, the proof that $F$ satisfies excision would end here,
but we do not.  The internal function space in $\bold{Space}$ is $k\Cal
C(\cdot,\cdot)$, where $k$ is the $k$-ification functor from the
category of topological spaces to $\bold{Space}$ and $\Cal
C(\cdot,\cdot)$ carries the compact-open topology.  By adjunction (at
the level of simplicial categories), $H=sk\Cal C(|\cdot|,Y)$ \cite{12,
II.3.14}, but $sk=s$, so $H=s\Cal C(|\cdot|,Y)$.  If $A\to B$ is a weak
equivalence of simplicial sets, then the weak equivalence $|A|\to|B|$
has a homotopy inverse, which induces a homotopy inverse for $k\Cal
C(|B|,Y)\to k\Cal C(|A|,Y)$, so $H(B)\to H(A)$ is a weak equivalence. 
Finally, since $H$ has a left adjoint \cite{12, II.2} (this is one of
the defining properties of a simplicial category), $H$ preserves limits,
i.e., takes colimits in $s\bold{Set}$ to limits in $s\bold{Set}$. 

Since $H$ preserves limits and weak equivalences and turns cofibrations
into fibrations, $H$ turns homotopy pushouts in $s\bold{Set}$, i.e.,
homotopy pullbacks in $s\bold{Set}^{\text{op}}$, into homotopy
pullbacks.  Also, the singular functor $s$ turns homotopy pushouts of
cofibrant spaces into homotopy pushouts in $s\bold{Set}$.  To prove that
$F$ satisfies excision, it therefore suffices to show that if $U$ and
$V$ are open subsets of $X$, then the square
$$\CD
U\cap V @>>> U \\
@VVV      @VVV \\
V  @>>>  U\cup V
\endCD $$
is a homotopy pushout.  To calculate the homotopy pushout of $U\gets
U\cap V\to V$, we factor $U\cap V\hookrightarrow U$ through its mapping
cylinder $M$ into a cofibration followed by a weak equivalence
\cite{22, 6.3}, and take the ordinary pushout of the diagram $M\gets
U\cap V\to V$.  The homotopy pushout $P$ turns out to be the product
$(U\cap V)\times [0,1]$ with $U$ glued to $(U\cap V)\times\{0\}$, and $V$
to $(U\cap V)\times\{1\}$.  We need to verify that the projection 
$P\to U\cup V$ is a weak equivalence.  It is easy to see that a section
of the projection is a homotopy inverse for it, and finding a section is
tantamount to finding a continuous function $U\cup V\to[0,1]$ equal to $0$
on $U\setminus V$ and $1$ on $V\setminus U$.  Since $U\cup V$ is normal,
such a function is provided by the Urysohn lemma.

Finally, since every open subset $U$ of $X$ is cofibrant, the weak
equivalence $|sU|\to U$ has a homotopy inverse, so the induced map
$s\Cal C(U,Y) \to s\Cal C(|sU|,Y)=F(U)$ is a weak equivalence.  Since
$F$ satisfies excision, so does $s\Cal C(\cdot,Y)$. 
\qed\enddemo

\specialhead 4. Complex manifolds as simplicial sheaves on the Stein
site  \endspecialhead

\noindent
Let $\Cal M$ be the category of complex manifolds (second countable but
not necessarily connected) and holomorphic maps.  In this section, we
shall embed $\Cal M$ into a model category, suitable for a
homotopy-theoretic interpretation of the Oka-Grauert property. 

Let $\Cal S$ be the category of Stein manifolds and holomorphic maps. 
This is a small category (or at least equivalent to one), since a
connected Stein manifold can be embedded into Euclidean space.  We view
$\Cal S$ as a site with the \lq\lq usual\rq\rq\  topology, in which a
cover of a Stein manifold $S$ consists of a family of isomorphisms onto
Stein open subsets of $S$ which cover $S$.  We only need to verify that
covers can be pulled back to covers \cite{21, III.2}.  This is implied
by the following lemma. 

\proclaim{4.1. Lemma} Let $f:X\to Y$ be a holomorphic map between
complex manifolds.  If $X$ is Stein and $V$ is a Stein open subset of
$Y$, then the preimage $f^{-1}(V)$ is Stein.  
\endproclaim

Note that the lemma shows that the intersection of finitely many Stein
open subsets of a complex manifold is Stein.  (The ambient manifold need
not be Stein.) Namely, for two subsets $U$ and $V$, the intersection is
the preimage of $V$ under the inclusion $U\hookrightarrow U\cup V$.  For
more sets, iterate this. 

\demo{Proof} Holomorphic functions separate points on $U=f^{-1}(V)$
since they do on $X$.  To show that $U$ is holomorphically convex, we
follow \cite{15, 2.5.14}.  Let $K\subset U$ be compact.  We
need to show that the holomorphic hull $\hat K_U$ of $K$ in $U$ is
compact.  Since $\hat K_X$ is compact and contains $\hat K_U$, it
suffices to show that $\hat K_U$ is closed in $U$.  Now $f(K)\subset V$
is compact, so $\widehat{f(K)}_V$ is compact.  If $h\in\Cal O(V)$ and
$x\in\hat K_U$, then $|h(f(x))| \leq \|h\circ f\|_K = \|h\|_{f(K)}$,
so $f(x)\in\widehat{f(K)}_V$.  Hence, $f$ maps the closure of $\hat K_U$
into $\widehat{f(K)}_V\subset V$, so the closure is in $U$.
\qed\enddemo

We shall refer to this topology on $\Cal S$ as the {\it fine} topology
and to the trivial topology, in which a cover consists of a single
isomorphism, as the {\it coarse} topology.  Note that a point, denoted
$\frak p$, is the final object in $\Cal S$, and the empty manifold
$\emptyset$ is the initial object.

\smallskip
Let $p$ be a point in a Stein manifold $S$.  If $F$ is a presheaf on
$\Cal S$, we define the {\it stalk} $F_p$ of $F$ at $p$ to be the
filtered colimit of the sets of sections $F(U)$, where $U$ is a Stein
neighbourhood of $p$ in $S$.  By restricting to the small site of each
Stein manifold, we see that the family of stalk functors $\cdot_p:\Shv\Cal
S\to\bold{Set}$, $F\mapsto F_p$, $p\in S\in\Cal S$, is faithful, meaning
that maps of sheaves are equal if they induce the same maps on all
stalks.  Since homotopy groups respect filtered colimits of simplicial
sets, we see that a map $f:F\to G$ of simplicial presheaves on $\Cal S$
is a fine weak equivalence (in the sense of Jardine) if and only if the
induced map of stalks $f_p:F_p\to G_p$ is a weak equivalence for all
$p\in S\in\Cal S$. 

Let us remark that the stalk functors just defined really are stalks (or
points of $\Cal S$) in the sense of topos theory: they have right
adjoints and preserve finite limits as functors $\Shv\Cal S\to\bold{Set}$
(so in particular, $\Cal S$ has enough points).  First of all, since
limits of sheaves are taken pointwise and finite limits commute with
filtered colimits in $\bold{Set}$ \cite{20, IX.2}, the stalk functor
$\cdot_p$ preserves finite limits.  It also preserves colimits.  A
colimit of sheaves is the sheafification of the pointwise colimit.  In
$\Cal S$, sheafification commutes with the restriction to the small site
of any Stein manifold, because all covers can be realized in the
manifold.  Hence, preservation of colimits can be reduced to the case of
a single manifold, where it holds by the standard result on topological
spaces (the Stein open subsets form a basis for the usual topology).  It
is easily verified that the functor $\Cal S\to\bold{Set}$ obtained by
restricting $\cdot_p$ to representable sheaves is both filtering and
continuous, so its Kan extension $\Shv\Cal S\to\bold{Set}$ is a point in
$\Cal S$ \cite{21, VII.5,6}, but since $\cdot_p$ preserves colimits, it
is its own Kan extension.  The right adjoint of $\cdot_p$ can in fact be
described explicitly \cite{21, VII.5}: it takes a set $A$ to the \lq\lq
skyscraper sheaf\rq\rq\  $X\mapsto \hom_{\bold{Set}}(\Cal
O(\cdot,X)_p, A)$ on $\Cal S$ (which does not really look like a
skyscraper at all).

\smallskip
We will embed $\Cal M$ into $s\Shv\Cal S$.  First note that by the
Yoneda lemma \cite{20, III.2}, there is a full embedding of $\Cal M$
into the category $\Pre\Cal M$ of presheaves of sets on $\Cal M$, given
by $X\mapsto\Cal O(\cdot, X)$.  It is easy to prove directly that this
is still true for the smaller site $\Cal S$. 

\proclaim{4.2. Proposition}  The Yoneda functor $\Cal M\to\Pre\Cal S$
is a full embedding.
\endproclaim

\demo{Proof} We need to show that the functor is both faithful and full,
i.e., that it induces bijections on all sets of morphisms.  Let
$f,g:X\to Y$ be maps in $\Cal M$ such that $f_*=g_*:\Cal O(\cdot,
X)\to\Cal O(\cdot, Y)$.  Let $\iota$ be the map $\frak p \to X$ with
image $\{ p\}$ for $p\in X$.  Since $f\circ\iota=g\circ\iota$, we have
$f(p)=g(p)$ and $f=g$. 

As for fullness, let $\alpha:\Cal O(\cdot, X)\to \Cal O(\cdot, Y)$ be a
map in $\Pre\Cal S$.  Now $\alpha:\Cal O(\frak p,X)\to\Cal O(\frak p,Y)$
gives a map $f:X\to Y$.  For a Stein manifold $S$ and $p\in S$ we have a
diagram
$$ \CD
\Cal O(S,X) @>\alpha>> \Cal O(S,Y) \\
@VVV                   @VVV  \\
\Cal O(\frak p, X) @>\alpha>> \Cal O(\frak p,Y)
\endCD $$
where the vertical arrows are induced by the map $\frak p\to S$ with
image $\{ p\}$.  Hence, for $h\in\Cal O(S,X)$, we have
$\alpha(h)(p)=f(h(p))$, so $\alpha=f_*$.  Finally, $f$ is holomorphic
because it preserves holomorphic maps from balls, mapping $\Cal O(B, X)$ 
into $\Cal O(B,Y)$, where $B$ is any open ball in Euclidean space.
\qed\enddemo

The Yoneda embedding restricts to a full embedding of $\Cal M$ into the
category of presheaves of topological spaces on $\Cal S$, if we equip
each set $\Cal O(S,X)$ with the compact-open topology.  Finally, let us
postcompose this functor with the singular functor.  This yields an
embedding $\Cal M\to s\Shv\Cal S$, taking a complex manifold $X$ to the
simplicial presheaf $s\Cal O(\cdot, X)$, which is clearly a sheaf with
respect to the fine topology on $\Cal S$ (recall that the singular
functor preserves limits).  We will not address the issue of fullness 
here, but only remark that every morphism $s\Cal O(\cdot, X)\to 
s\Cal O(\cdot, Y)$ is clearly given by a holomorphic map $X\to Y$ at 
the level of vertices.  We will often view a complex manifold $X$ as an
object of $s\Shv\Cal S$ and write $X$ for $s\Cal O(\cdot, X)$. 

If $p$ is a point in a Stein manifold $S$ and $\dim_p S=m\geq 0$, then
the stalk of a complex manifold $X$ at $p$ is simply the colimit as
$n\to\infty$ of $s\Cal O(\frac 1 n\Bbb B_m,X)$, where $\Bbb B_m$ is the
open unit ball in $\Bbb C^m$ and the maps between the scaled balls
$\frac 1 n\Bbb B_m$ are the inclusions.  This is in fact a homotopy
colimit, since (by the identity theorem) all the maps are cofibrations
(just dualize the theory of towers in \cite{12, VI.1}).  Now for any
$m\geq 0$ and $r>0$, $\Cal O(r\Bbb B_m,X)$ is weakly equivalent to $X$
itself, so all the stalks of $X$ are weakly equivalent to $sX$. 

A holomorphic map $X\to Y$ is a cofibration (with respect to either
topology) if and only if it is injective.  It is a coarse weak
equivalence if and only if it induces a topological weak equivalence
$\Cal O(S,X)\to\Cal O(S,Y)$ for every Stein manifold $S$.  It is a fine,
i.e., stalkwise, weak equivalence if and only if it is a topological
weak equivalence.  As for fibrations, we only remind the reader that 
they are defined by a right lifting property with respect to acyclic 
cofibrations. 

\smallskip\noindent
{\bf Examples.}  Since $\Bbb D$ and $\Bbb C$ are holomorphically
contractible, they are both coarsely weakly equivalent to a point.  The
same holds for any star-shaped domain in Euclidean space.
The inclusion $\Bbb D^\times\hookrightarrow\Bbb C^\times$ is a fine but
not a coarse weak equivalence.  Namely, by Liouville's theorem, $\Cal
O(\Bbb C^\times, \Bbb D^\times)=\Bbb D^\times$, but $\Cal O(\Bbb
C^\times, \Bbb C^\times)$ has infinitely many connected components, one
for each winding number about the origin.

\smallskip
Now let $X$ be a complex manifold.  The inclusion $s\Cal
O(\cdot,X)\hookrightarrow s\Cal C(\cdot,X)$ is a fine weak equivalence
of simplicial sheaves on $\Cal S$ because every cover has a refinement
consisting of sets at which the inclusion is a weak equivalence: take a
refinement by balls, for instance.  The Oka-Grauert property is
satisfied by $X$ if and only if this inclusion is a coarse weak
equivalence.  Theorem 2.1 states that this is equivalent to the
simplicial sheaf $s\Cal O(\cdot,X)$ satisfying {\it finite excision} in
the sense of Section 2, meaning that for every finite cover
$\{U_1,\dots,U_m\}$ in $\Cal S$, $s\Cal O(\bigcup U_i,X)$ is not only
the limit but also the homotopy limit of the diagram whose objects are
the simplicial sets $s\Cal O(U_{i_1}\cap\dots\cap U_{i_k}, X)$ for $1
\leq i_1 < \dots < i_k \leq m$ and whose arrows are induced by
restriction maps.  This property may also be expressed by saying that
$s\Cal O(\cdot,X)$ is a {\it finite homotopy sheaf}.  We can now state
Gromov's Oka principle and our interpretation of its conclusion, the
Oka-Grauert property, as follows. 

\proclaim{4.3. Theorem}  Let $X$ be a complex manifold.  
\roster
\item The fine weak equivalence $s\Cal O(\cdot,X)\hookrightarrow s\Cal
C(\cdot,X)$ of simplicial sheaves on the Stein site is coarsely acyclic 
if and only if $s\Cal O(\cdot,X)$ is a finite homotopy sheaf.  
\item If $X$ has a spray, then $X$ represents a finite homotopy sheaf on
the Stein site.
\endroster
\endproclaim

Finite excision for pointwise fibrant simplicial presheaves on $\Cal S$
is clearly invariant under coarse weak equivalences.  Let us show that
descent implies finite excision.  As remarked before Theorem 3.1,
descent implies two-set excision, but since we are not working on a
topological space (the union of Stein open subsets is usually not
Stein), we cannot simply refer to Theorem 3.3 to get finite excision.  A
simplicial presheaf satisfying descent is coarsely weakly equivalent to
a finely fibrant simplicial sheaf $F$ on $\Cal S$, so it suffices to
show that $F$ satisfies finite excision.  Let $S$ be a Stein manifold. 
Via its \ac etale space, $F|S$ extends to a simplicial sheaf $\tilde F$
on $S$ with its usual topology.  It suffices to show that $\tilde F$ is
flabby; we then invoke Theorem 3.3.  By the Yoneda lemma, if $X$ is a
Stein manifold, then $F(X)=\bold{Hom}_{s\Shv\Cal S}(\hat\Cal O(\cdot,
X), F)$, where the sheaf $\Cal O(\cdot,X)=\hom_{\Cal S}(\cdot,X)$ of sets
on $\Cal S$ has been turned into a simplicial sheaf $\hat\Cal
O(\cdot,X)$ in the trivial way (the same set in all degrees; all face
and degeneracy maps are the identity).  If $V$ is an open subset of $S$,
let $\Cal B$ be a Stein basis for the topology of $V$, viewed as a
subdiagram of the site of $S$.  Then $\Cal O(\cdot, V)$ is the sheaf
colimit of the diagram $\Cal O(\cdot, \Cal B)$, and
$$\bold{Hom}_{s\Shv\Cal S}(\hat\Cal O(\cdot, V), F)
= \lim\bold{Hom}_{s\Shv\Cal S}(\hat\Cal O(\cdot, \Cal B), F)
= \lim F(\Cal B)=\tilde F(V).$$
If $W\subset V$ are open subsets of $S$, then the induced map $\hat\Cal
O(\cdot,W)\to\hat\Cal O(\cdot,V)$ is clearly a cofibration (pointwise 
injection), so by Quillen's Axiom SM7, the restriction map 
$\tilde F(V)\to\tilde F(W)$ is a fibration, and the proof is complete.

Note, finally, that a nondiscrete complex manifold $X$ is never flabby,
let alone coarsely or finely fibrant.  Namely, let $D_r$ be the open
disc of radius $r$ centred at the origin in the complex plane.  The
inclusion $D_1\hookrightarrow D_2$ is a monomorphism in $\Cal S$, but
since there are holomorphic maps $D_1\to X$ that do not extend
holomorphically to $D_2$, it is easily seen that the restriction map
$s\Cal O(D_2,X)\to s\Cal O(D_1,X)$ is not a fibration.

\specialhead 5. Partial descent on the quasi-projective site 
\endspecialhead

\noindent It is natural to ask whether a finite homotopy sheaf on $\Cal
S$ satisfies descent.  This would turn Gromov's Oka principle into a
descent theorem, somewhat analogous to such results as Brown-Gersten
descent and (unstable) Nisnevich descent in algebraic geometry.  I do
not know the answer: the finiteness properties that make descent
possible in algebra --- the Zariski topology being binoetherian,
essentially --- do not hold in analysis.  We hope to address this question
in future work.  In the meantime, let us show how Brown-Gersten descent
easily implies partial descent of sorts for quasi-projective manifolds 
with the Oka-Grauert property. 

Let $\Cal A$ be the category of quasi-projective complex manifolds,
i.e., smooth Zariski open sets in projective varieties, and algebraic
maps.  We put the usual Zariski topology on $\Cal A$ by defining a cover
of a quasi-projective manifold $X$ to be a family of isomorphisms onto
Zariski open subsets of $X$ which cover $X$ (covers can be pulled back
to covers because algebraic maps are Zariski continuous).  Taking $X$
to the sheaf of sets $\Cal O(\cdot,X)$ on the small site $\Cal A$
defines an embedding of $\Cal A$ into $\Shv\Cal A$ (it is faithful
because the Yoneda embedding $X\mapsto\hom_{\Cal A}(\cdot, X)$ is).  As
before, we equip each set of holomorphic maps with the compact-open
topology, apply the singular functor, and obtain an embedding of $\Cal
A$ into the model category $s\Shv\Cal A$ of simplicial sheaves on $\Cal
A$, taking a quasi-projective manifold $X$ to the simplicial sheaf
$s\Cal O(\cdot, X)$.  This embedding is not full, but every morphism
$s\Cal O(\cdot, X)\to s\Cal O(\cdot,Y)$ is given by a holomorphic map
$X\to Y$ at the level of vertices. 

Now let $X$ be a quasi-projective manifold.  We claim that the
simplicial sheaf $G=s\Cal C(\cdot, X)$ on $\Cal A$ satisfies descent. 
Let $G\to F$ be a fibrant model for $G$, i.e., an acyclic cofibration to
a fibrant simplicial sheaf $F$ on $\Cal A$.  Let $A$ be a
quasi-projective manifold with the Zariski topology.  By Theorem 3.4,
$G|A$ satisfies excision.  So does $F|A$, since it is flabby.  Hence,
the weak equivalence $G\to F$ is pointwise acyclic by Proposition 3.2,
and $G$ satisfies descent. 

Suppose now that $X$ has the Oka-Grauert property, so $s\Cal
O(S,X)\hookrightarrow G(S)$ is a weak equivalence for every Stein
manifold $S$.  Every cover in $\Cal A$ has a refinement consisting of
Stein Zariski open sets, so $s\Cal O(\cdot,X)\to G$ is an acyclic
cofibration between simplicial sheaves on $\Cal A$.  Hence, the
composition $s\Cal O(\cdot, X)\to F$ is a fibrant model for $s\Cal
O(\cdot, X)$, and $s\Cal O(S,X)\to F(S)$ is a weak equivalence for every
Stein manifold $S$ in $\Cal A$.  Since any two fibrant models for the
same object in $s\Shv\Cal A$ are pointwise weakly equivalent, this holds
for every fibrant model for $s\Cal O(\cdot, X)$, and we have proved the
following \lq\lq Stein descent theorem\rq\rq. 

\proclaim{5.1. Theorem}  Let $X$ be a quasi-projective manifold and $F$
be a fibrant model for $s\Cal O(\cdot,X)$ in $s\Shv\Cal A$.  If $X$ has
the Oka-Grauert property, then $s\Cal O(S,X)\to F(S)$ is a weak
equivalence for every quasi-projective Stein manifold $S$.
\endproclaim

\enddocument